\documentclass[a4paper, 10pt, ]{article}

\usepackage{amssymb,amsmath,enumerate,graphicx,hyperref,color}
 
\usepackage{t1enc}
 
\usepackage[T1]{fontenc}

\usepackage{amsfonts}

\usepackage{color}
\usepackage{graphicx}
\usepackage{rotating}
\usepackage{caption}
\usepackage{fancyhdr}
 
\usepackage{hyperref}
 
\newtheorem{lemma}{Lemma}[section]
\newtheorem{proposition}[lemma]{Proposition}
\newtheorem{theorem}[lemma]{Theorem}

\newtheorem{remark}[lemma]{Remark}

\newtheorem{example}[lemma]{Example}

\newtheorem{assumption}[lemma]{Assumption}

\newcommand{\beq}{\begin{eqnarray}}
\newcommand{\enq}{\end{eqnarray}}
\newcommand{\be}{\begin{eqnarray*}}
\newcommand{\en}{\end{eqnarray*}}
\newcommand{\ben}{\begin{eqnarray*}}
\newcommand{\enn}{\end{eqnarray*}}

\newcommand{\R}{\mathbb R}

\newcommand{\E}{\mathbb E}

\newcommand{\1}{{\mathbf I}}

\newcommand{\F}{\mathbb F}

\let\cal=\mathcal

\newcommand{\Ac}{\mathcal{H}} 
\newcommand{\Fc}{\mathcal{F}}

\newcommand{\eps}{\varepsilon}  
  
\newcommand{\x}{\times}

\def \vp{\varphi}
\def\Q{{\mathbb Q}}

 \def\xr{\omega}

  \def\E{{\mathbb E}}
    \def\P{{\mathbb P}}

	\def\Om{\Omega}
	\def\om{\omega}

	\def\0{\mathbf{0}}

 \def\Mc{{\cal M}}

  \def\Pc{{\cal P}}

  \def\R{{\mathbb R}}

 \def\Cb{{\mathbb C}_{\rm r}}

  \def\Yb{\overline Y}

\title{A quasi-sure optional decomposition and super-hedging result on the Skorokhod space}
\author{
Bruno Bouchard
\footnote{CEREMADE, Universit\'e Paris-Dauphine, PSL, CNRS.  bouchard@ceremade.dauphine.fr. } 
\and 
Xiaolu Tan
\footnote{Department of Mathematics, The Chinese University of Hong Kong. xiaolu.tan@cuhk.edu.hk.}} 

\date{\today}

\begin{document}

\maketitle

\begin{abstract} 
	We prove a robust super-hedging duality result for path-dependent options on assets with jumps, in a continuous time setting. It requires that the collection of martingale measures is rich enough and that the payoff function satisfies some continuity property. It is a by-product of a quasi-sure version of the optional decomposition theorem, which can also be viewed as a functional version of It\^{o}'s Lemma, 
	that applies to non-smooth functionals (of c\`adl\`ag processes) which are only concave in space and non-increasing in time, in the sense of Dupire. 
\end{abstract}

\section{Introduction}

	A key element in the proof of the  super-hedging duality is the optional decomposition theorem.
	Let $X$ be a stochastic process on some probability space and
	 consider the class of all equivalent martingale measures under which $X$ is a local martingale.
	Let $V$ be a  supermartingale under all these equivalent martingale measures. Then, the classical optional decomposition theorem states that there exists a predictable process $H$ and a non-decreasing process $C$
	such that $V = V_0 + \int_0^{\cdot} H_r\cdot dX_r - C$, almost surely.
	Initially introduced by El Karoui and Quenez \cite{ELQ} in the case where $X$ has continuous paths,
	it was then extended to the c\`adl\`ag paths case in Kramkov \cite{Kramkov}, F\"ollmer and Kabanov \cite{FK97}, Delbaen and Schachermayer \cite{DSch}, F\"ollmer and Kramkov \cite{FollKramkov}.

	\vspace{0.5em}

	The optional decomposition theorem has been recently studied in the robust context in which negligible sets, associated to one reference measure, are replace by polar sets associated to  a family $\Pc$ of (singular) reference measures.
	While the classical optional decomposition theorem ensures the existence of a couple $(H^{\P},C^{\P})$ for each $\P \in \Pc$, the robust version consists essentially in aggregating the family $(H^{\P})_{\P \in \Pc}$ into a universal process $H$, independent of $\P$.
	As in the classical case, the robust optional decomposition theorem is key to prove (and is motivated by) the robust super-hedging duality.
	For the literature, let us refer to Bouchard and Nutz \cite{BN13} for  discrete time models,  to Biagini, Bouchard, Kardaras and Nutz \cite{BBKN13},
	Neufeld and Nutz \cite{NN13}, Possama\"i, Royer and Touzi \cite{PRT} for the continuous time models when $X$ has continuous paths, and to 
	Nutz \cite{nutz2015OptionalDecompo} when $X$ has  c\`adl\`ag paths, among others.	
	
	\vspace{0.5em}
	
	In this paper, we aim at providing a robust optional decomposition theorem
	when  $X$ is a c\`adl\`ag process, and the super-martingale $V$ is given as a concave functional of $X$ which decreases in time (both in the sense of Dupire \cite{dupireito}, see below). As mentioned above, the case where $X$ has c\`adl\`ag paths has already been studied in Nutz \cite{nutz2015OptionalDecompo}. However, it uses the crucial condition   that the jump part of $X$ is ``dominated'' by its diffusion part. Basically, this allows to write the decomposition under each $\P$ and then characterize $H^{\P}$ in terms of the quadratic variations of the continuous parts of $V$ and $X$. Because the latter can be defined pathwise, this allows one to show that all the $H^{\P}$ admit a common version, universally defined.

	\vspace{0.5em}
	
	We do not impose such a domination condition, but we instead require some continuity, concavity and monotonicity  on $(t,\omega)\mapsto V(t,\omega)$. 
	It is enough to derive a  version of It\^{o}'s formula\footnote{It could also be viewed as a version of Meyer-Tanaka's formula, except that  the bounded variation part is not identified in terms of the local time processes.} solely expressed in terms of the first order horizontal Dupire's derivative (or elements of the associated super-differential), or, said differently, a robust optional decomposition.
	In the robust super-hedging problem, the supermartingale $V$ is obtained as the sup over a family of martingale measures of the expectation of the  payoff,
	and the above conditions are satisfied as soon as the family of martingale measures is rich enough, and the payoff function enjoys some continuity. It thus suffices to apply this decomposition to deduce the super-hedging duality, and that a super-hedging strategy is actually associated to the super-hedging price. 
	
	\vspace{0.5em}

As a by-product, we prove that any locally-bounded path-dependent Dupire-concave function of a $\R^{d}$-valued semi-martingale remains a semi-martingale, thus generalizing the result of  Meyer \cite[Chapter VI]{Meyer} (see also Carlen and Protter \cite{CarlenProtter}).  
	
	\vspace{0.5em}

	The rest of this paper is organized as follow.  We first introduce some notations that will be used all over this paper. We state our version of the robust optional decomposition theorem in Section \ref{sec:OptionalDecomp}. In Section \ref{sec:Duality}, we  provide a pretty general version of the robust super-hedging duality for continuous payoffs, including a typical example in which the components of $X$ are restricted to remain non-negative.

	\vspace{0.5em}

	{\bf Notations.} 
	$\mathrm{(i)}$. 	
	Let $E \subseteq \R^d$ be a closed convex set,
	we denote by $\Om = D([0,T], E)$ be the space of all c\`adl\`ag $E$--valued paths on $[0,T]$, with canonical filtration $\F = (\Fc_t)_{0 \le t \le T}$ and canonical process $X(\om) :=\om$. 
	We endow $\Om$ with the sup-norm topology induced by $\|\xr-\xr'\|:=\sup_{t\in [0,T]}|\xr_{t}-\xr'_{t}|$ for $\xr,\xr'\in \Om$. 
	For $(t, \xr)\in \Theta:= [0,T] \x \Om$,
	we consider the (optional) stopped path $\xr_{t \wedge \cdot}:=(\xr_{t\wedge s})_{s\in [0,T]}$, 
	and (predictable) stopped path $\xr^{t-} := (\xr^{t-}_s)_{s \in [0,T]}$ defined by $\xr^{t-}_s := \xr_s \mathbf{1}_{\{ s \in [0,t) \}} + \xr_{t-} \mathbf{1}_{\{s \in [t,T]\}}$.
	A function $\vp: \Theta \to \R$ is said to be non-anticipative if $\vp(t, \xr)=\vp(t, \xr_{t \wedge \cdot})$ for all $(t,\xr)\in \Theta$.
	
	\vspace{0.5em}

	$\mathrm{(ii)}$. For a function $\vp: \Theta \to \R$, we follow Dupire \cite{dupireito} (see also Cont and Fourni\'e \cite{cont2013functional}) to introduce the Dupire derivatives as follows:
	$\vp$ is said to be horizontally differentiable if,
	for all $(t,\xr)\in [0,T) \x \Om $, its horizontal derivative
	$$
		\partial_{t}\vp(t,\xr) 
		~:=~ 
		\lim_{h\searrow 0} \frac{\vp(t+h,\xr_{t \wedge \cdot}) - \vp(t,\xr_{t \wedge \cdot})}{h}
	$$
	is well-defined,
	$\vp$ is said to be vertically differentiable if, for all $(t,\xr)\in \Theta$, the function
	$$
		y \in E \longmapsto \vp(t,\xr\oplus_{t} y) \in \R
		~\mbox{is differentiable, with}~
		\xr\oplus_{t} y := \xr \1_{{[0,t)}}+\1_{[t,T]} y,
	$$
	whose derivative at $y= \om_t$ is defined as the vertical derivative $\nabla_{\xr} \vp(t,\xr)$ of $\vp$. 

	\vspace{0.5em}

	$\mathrm{(iii)}$. 
	Let us denote by $\Cb(\Theta)$ the class of all non-anticipative functions $\vp: \Theta \longrightarrow \R$ such that
	$$
		\vp(t^{n},\xr^{n}) \longrightarrow \vp(t,\xr)
		~~\mbox{whenever}~~
		t^{n} \ge t,
		~t^n \longrightarrow t
		~\mbox{and}~
		\big\| \xr^n_{t^n \wedge \cdot} - \xr_{t \wedge \cdot} \big\| \longrightarrow 0.
	$$
	We say that $\vp\in \Cb^{0,1}(\Theta)$ if both $\vp$ and $\nabla_{\xr}\vp$ are well defined and belong to $\Cb(\Theta)$.

	\vspace{0.5em}

	$\mathrm{(iv)}$.
	A non-anticipative map $\vp : \Theta \longrightarrow \R$ is said to be Dupire-concave if,
	for all $t \in [0,T]$,  $\om^1, \om^2 \in \Om$, such that $\om^1 = \om^2$ on $[0,t)$, and $\theta \in [0,1]$,
	\begin{equation} \label{eq:def_concave}
		\vp(t,  \theta \om^1+(1-\theta) \om^2) ~\ge~ \theta \vp(t, \om^1) + (1-\theta) \vp(t, \om^2).
	\end{equation}
	For a Dupire-concave function $\vp$, one can define the Dupire super-differential (set)
	$$
		\partial \vp(t, \xr)
		~:=~
		\big \{z\in \R^{d}: \vp(t,\xr\oplus_{t} y)\le \vp(t,\xr)+z\cdot (y - \om_t),\;\forall\; y \in E \big\}. 
	$$
	The map $\vp$ is said to be Dupire-non-increasing in time if 
	$$
		\vp(t, \xr_{t \wedge \cdot}) \ge \vp(t+h,  \xr_{t \wedge \cdot}),\;\mbox{ for all $(t, \xr) \in \Theta$ and $h \in [0, T-t]$.}
	$$
	Note that a map $\vp : \Om \to \R$ can be associated to the Dupire-non-increasing map $t\mapsto \vp(\om_{t\wedge \cdot})$.

	\vspace{0.5em}

	\noindent$\mathrm{(v)}$. 
	Given a locally bounded predictable process $H$, and a (c\`adl\`ag) semimartingale $X$, we write (as usual) $\int_s^t H_{r}\cdot dX_{r}$ for the stochastic integral $\int_{(s,t]} H_{r} \cdot dX_{r}$.
	In the case that the law of $X$ depends on the reference probability measure $\P$, the integral $\int_s^t H_{r}\cdot dX_{r}$ depends also on $\P$, 
	which is usually omitted whenever it is obviously given by the context.

\section{Optional decomposition of Dupire concave functionals}	
\label{sec:OptionalDecomp}

	We provide immediately our version of the optional decomposition theorem, which is the key ingredient for proving the super-hedging duality of Theorem \ref{thm:surrep duality} below. 
	It can also be seen as a functional version of It\^{o}'s or Meyer-Tanaka's formula, as it generalizes both up to the fact that the bounded variation part entering our decomposition is not explicitly characterized.
	As opposed to the classical versions of the optional decomposition theorem mentioned in the introduction, it is a functional one as our starting point is not that $t\mapsto V(t,\omega)$ is a super-martingale   under martingale measures (although one can easily check that our assumptions imply this). 
 
	\vspace{0.5em}

	Recall that $\Om := D([0,T], E)$ denotes the canonical space of all $E$-value  c\`adl\`ad paths on $[0,T]$.
	Throughout the paper, we assume that $E \subseteq \R^d$ is closed convex set, {with non-empty interior},
	and moreover that there exists a compactly supported smooth density function $\phi: \R^d \to \R$ such that
		\begin{equation} \label{eq:Support_phi}
			\mbox{the map}~
			y \longmapsto \phi(y-x)
			~\mbox{is supported in}~E,
			~\mbox{for all}~
			x \in E.
		\end{equation}
	Let $\Pc$ denote the collection of all Borel probability measures on $\Om$, 
	under which the canonical process $X$ is a semimartingale,
	that is, a (c\`adl\`ag) process which can be decomposed as the sum of a local martingale and an adapted finite variation process, w.r.t.~the (augmented) canonical filtration.
	For $s \in [0,T]$, denote also 
	$X^{s-}_r(\om) := \om^{s-}_r$, or equivalently, $X^{s-}_r := X_{s \wedge r} - \Delta X_s \1_{\{r \ge s\}}$ for all $r \in [0,T]$.
	
	\begin{theorem} \label{thm:OptionalDecomp}
		Let $V \in \Cb(\Theta)$ be Dupire-concave, Dupire-non-increasing in time,
		and such that
		\begin{equation} \label{eq:bound_V}
			\sup \big\{
				|V(t, \om)| + |z|
				~: (t,\om) \in \Theta, ~\|\om\| \le K, ~z \in \partial V(t,\om)
			\big\}
			< \infty,
			~~\mbox{for all}~
			K > 0.
		\end{equation}
		Then, there exists a $\F$--predictable locally bounded process $H: \Theta \to \R^d$, together with
		a collection of non-decreasing processes $\{C^{\P} :\P\in \Pc\}$,
 		satisfying
		\begin{equation} \label{eq:OptDecomp}
			V(t, X)= V(0, X) + \int_{0}^t H_s \cdot d X_s - C^{\P}_t, ~t \in [0,T],  ~\P{\rm-a.s.} \;\forall\; \P\in \Pc.
		\end{equation}
		Moreover, $H_{s}  \in \partial V(s, X^{s-})$ for all $s\in [0,T]$, $\Pc$-q.s.
	\end{theorem}
	
	Let us make some remarks before proving this result. 
	
	\begin{remark} \label{rem: H explicit}
		An explicit formula for $H$ is given in Remark \ref{rem:DVplus} as an element of the super-differential of $\partial V(\cdot,X^{\cdot-})$.

		\vspace{0.5em}
		
		{\rm (i)} From this point of view, it can be considered as a  version of the functional Meyer-Tanaka's formula, except that $C^{\P}$ is not identified to be associated to  local time processes.
		In particular, when $E = \R^d$ and $V(t, \om) = f(\om_t)$ for some convex function $f: \R^d \to \R$, it satisfies clearly all the conditions in Theorem \ref{thm:OptionalDecomp},
		and the decomposition result \eqref{eq:OptDecomp} implies the result of Meyer \cite[Chapter VI]{Meyer} (see also Carlen and Protter \cite{CarlenProtter}) which states that a convex function of a semimartingale is still a semimartingale.  Our result provides a path-depend version of this (apply it to $V(t,\omega):=f(\omega_{t \wedge \cdot})$ with $f: \Om \to \R$). 
		In the functional (path-dependent) case, when $X$ is a one-dimensional process with continuous paths, 
		such a decomposition has been derived in Saporito \cite{Saporito} with an explicit expression of $C^{\P}$
		in terms of the local times of $X$, but under additional smoothness conditions
		(see also Bouchard and Tan \cite{BT19} for a version when $V$ is only in  ${\Cb^{0,1}}(\Theta)$).		
		\vspace{0.5em}
		
		{\rm (ii)} It is clear that one may have different versions of the process $H$, which depends on the kernel $\phi$ used in part (ii) of the proof of Theorem \ref{thm:OptionalDecomp}.
		
		\vspace{0.5em}
		
		{\rm (iii)} Such an explicit formula is not available in the approach of  Nutz \cite{nutz2015OptionalDecompo} because it is based on the aggregation argument mentioned in the introduction (and does not assume any continuity).
	\end{remark}

	\begin{remark}  
		Note that a similar decomposition holds if $V+R$ satisfies the conditions of Theorem \ref{thm:OptionalDecomp} for some function $R$ that is ${\mathbb C}^{1,2}$ in the sense of Dupire,
		 using the functional It\^o's formula in Cont and Fourni\'e \cite{cont2013functional}).
		The difference will be that $\{C^{\P}:\;\P\in \Pc\}$ might not be a family of non-decreasing processes, unless, for instance, each element $\P\in \Pc$ admits an equivalent local-martingale measure for $X$ and $t\mapsto V(t,X)$ is a super-martingale under each of this local-martingale measures. 
	\end{remark}

	\begin{remark} 
		In Theorem \ref{thm:OptionalDecomp}, the stochastic integrals $(\int_{0}^t H_s \cdot d X_s)_{t\le T}$ depend on the reference measure $\P$, the fundamental point being that $H$ does not. However, following Nutz \cite{nutz2012pathwise}, these stochastic integrals could be aggregated into a single $\F^{*}$-optional process, with $\F^{*}$ defined as the universally augmented filtration. In this case, the corresponding non-increasing processes $\{C^\P :\P\in\Pc\}$ can also be aggregated into a process independent of $\P$. But, this requires to work under the Zermelo-Fraenkel set theory with the axiom of choice (ZFC) plus the Continuum Hypothesis, 
		as well as to assume the existence of a uniform dominating measure for the characteristics of $X$ (see \cite[Assumption 2.1]{nutz2012pathwise}).
	\end{remark}

	\begin{remark}
		{\rm (i)} Examples of  sets $E$ satisfying \eqref{eq:Support_phi} could be $E = \R^d$, or $E = \R_+^d$ with $\R_+ = [0, \infty)$, {or any cone of $\R^{d}$ (with non-empty interior)}.
		\vspace{0.5em}
		
		 {{\rm (ii)} Also note that the condition \eqref{eq:Support_phi} is only used to regularize $V$ into a function with continuous first order vertical Dupire derivative. It is not necessary if $\partial V(t,\omega)$ admits  a unique element for all $(t,\omega)\in \Theta$. In this case, the vertical Dupire derivative inherits the regularity of $V$ automatically, and there is no need for the intermediate smoothing procedure in part (ii) of the proof of Theorem \ref{thm:OptionalDecomp}.}
		\vspace{0.5em}
		
		{\rm (iii)} When $E = \R^d$, one has
		\begin{equation*}
			\sup \big\{|z| ~:z\in \partial \vp(t, \xr) \big\}
			~\le~ 
			\sup_{|y|\le 1} \big|\vp(t,\xr\oplus_{t} y)-\vp(t,\xr) \big|,
		\end{equation*}
		so that Condition \eqref{eq:bound_V} is equivalent to assuming that $V: \Theta \to \R$ is a locally bounded function.

	\end{remark}

	\begin{proof} \!\!\! {\bf (of Theorem \ref{thm:OptionalDecomp})} 
	$\mathrm{(i)}$ 
		Let us first assume in addition that $V \in \Cb^{0,1}(\Theta)$,
		so that $\nabla_{\om} V \in \Cb(\Theta)$
		and, for each $(t, \xr) \in \Theta$, $\nabla_{\xr} V(t, \xr)$ is the unique element in  $\partial V(t,\xr)$,
		or equivalently, in the super-differential of the map $y \mapsto V(r, \xr \oplus_r y)$ at $y=\om_r$.

		\vspace{0.5em}

		{\rm (a)} Let us fix 
		$s<t$ and consider a sequence of deterministic  discrete time grids $(\pi_n)_{n \ge 1}$,
		where $\pi_n = \{ t^n_k \}_{0 \le k \le n}$ satisfies 
		$$
			s= t^n_0 < t^n_1 < \cdots < t^n_{n} = t
			~~\mbox{and}~~
			|\pi_n| := \max_{k=1, \cdots, n} ( t^n_k - t^n_{k-1}) \longrightarrow 0,
			~\mbox{as}~ n \longrightarrow \infty.
		$$
		Next, with fixed $\delta>0$, we define the sequences of $\F$--stopping times $(\tau^{n}_k)_{k \ge 1}$, $n\ge 1$,  by 
		$$
			\tau^{n}_0 \equiv  s,
			~~\mbox{and}~~
			\tau^{n}_{k+1} := \inf \big\{r  {>} \tau^{n}_k ~: |\Delta X_r | \ge \delta  {\mbox{ or } r\in \pi_{n}} \big\},
			~k \ge 0.
		$$
		Observe that  the random number $m_n:=\max\{k\ge 0 : \tau^{n}_{k}\le t\}$ is finite, but is not uniformly bounded in general.

		\vspace{0.5em}
		
		For each $n \ge 1$ and $u \in [s, t]$, let the processes $X^n$ and $X^{n, u-}$
		be defined by 
		$$
			X^n_r := 
			\sum_{k=0}^{ {m_{n}}-1} 
				X_{\tau^n_k} \mathbf{1}_{\{r \in [\tau^n_k, \tau^n_{k+1})\}}
			+
			X_t \mathbf{1}_{\{r=t \}},
			~~
			X^{n, u-}_r := X^n_{u\wedge r} - \Delta X_u \mathbf{1}_{\{r \ge u \}},
			~
			r \in [s,t],
		$$
		so that
		$$
			X^n_{\tau^n_{k}} = X_{\tau^n_{k}}
			~~\mbox{and}~~
			X^{n, \tau^n_{k+1}-}_{\tau^n_{k+1}} = X_{\tau^n_{k+1}-},
			~~\mbox{for each}~ k \ge 0.
		$$
		Recall that $V$ is non-anticipative, Dupire-concave and that $V(\tau^n_k ,X^n) \ge V(\tau^n_{k+1}, X^n_{\tau^n_k \wedge \cdot})$ since it is Dupire-non-increasing in time. 
		It follows that 
		\begin{equation} \label{eq:Diff_k_p}
			V(\tau^n_{k+1}, X^n) - V(\tau^n_k ,X^n) 
			~\le~
			\nabla_{\xr} V(\tau^n_{k+1}, X^n_{\tau^n_k \wedge \cdot}) \cdot (X_{\tau^n_{k+1}}  - X_{\tau^n_k}),
			~~\mbox{if}~ \tau^n_{k+1} \in \pi_n,
		\end{equation}
		and
		\begin{eqnarray} \label{eq:Diff_k}
			&&
			V(\tau^n_{k+1}, X^n) - V(\tau^n_k ,X^n) \nonumber \\
			&\le&
			V(\tau^n_{k+1}, X^n) - V(\tau^n_{k+1}, X^{n, \tau^n_{k+1}-})
			+ V(\tau^n_{k+1}, X^{n, \tau^n_{k+1}-}) - V(\tau^n_{k+1}, X^n_{\tau^n_k \wedge \cdot}) \nonumber \\
			&\le&
			\nabla_{\xr} V(\tau^n_{k+1}, X^{n, \tau^n_{k+1}-}) \cdot \Delta X_{\tau^n_{k+1}}
			+
			\nabla_{\xr} V(\tau^n_{k+1}, X^n_{\tau^n_k \wedge \cdot}) \cdot (X_{\tau^n_{k+1}-}  - X_{\tau^n_k}),
			~~\mbox{if}~ \tau^n_{k+1} \notin \pi_n.~~~
		\end{eqnarray}	
		By summing up the two sides of \eqref{eq:Diff_k_p} and \eqref{eq:Diff_k} for $k=0, \cdots, m_n-1$, it follows that
		\begin{equation} \label{eq:OptDecomp_Xn}
			V(t, X^n) -  V(s, X^n)  
			\le
			I^{\delta}_n 
		\end{equation}
		with
		\begin{eqnarray*}
			I^{\delta}_n 
			\!\!\! &:=& \!\!\!
			\sum_{k=0}^{m_n-1} \Big(
				\nabla_{\xr} V(\tau^n_{k+1}, X^n_{\tau^n_k \wedge \cdot}) \cdot (X_{\tau^n_{k+1}}  - X_{\tau^n_k})
			\Big)  \1_{\{\tau^n_{k+1} \in \pi_n\}}
			\\
			&&\!\!\!\!\!\!\!\!\!\!\!\!
			+ 
			\sum_{k=0}^{m_n-1} \!\! \Big(
				\nabla_{\xr} V(\tau^n_{k+1}, X^{n, \tau^n_{k+1}-}) \cdot \Delta X_{\tau^n_{k+1}}
				+
				\nabla_{\xr} V(\tau^n_{k+1}, X^n_{\tau^n_k \wedge \cdot}) \cdot (X_{\tau^n_{k+1}-}  - X_{\tau^n_k})
			\Big) \1_{\{\tau^n_{k+1} \notin \pi_n\}} ,~~
		\end{eqnarray*}
		where we add the superscript $\delta$ on $I^{\delta}_n$ to emphasis the dependence of the random number $m_n$ and the stopping times $(\tau^n_k)_{1 \le k \le m_n}$ on $\delta > 0$.
		The term $I^{\delta}_n$ can be written as an integral w.r.t. $X$, but the integrand may not be adapted to the filtration $\F$.
		This motivates us to introduce 
		$$
			I_n 
			~:=~
			\sum_{k=0}^{n-1} \nabla_{\xr} V(t^n_{k+1}, X^n_{t^n_k \wedge \cdot}) \cdot \big( X_{t^n_{k+1}} - X_{t^n_k} \big)
			~=~
			\int_s^t H^n_r dX_r,
		$$
		where $H^n$ is the $\F$-predictable process defined by
		\begin{align*}
			H^n 
			:= 
			\sum_{k=0}^{n-1} 
			\nabla_{\xr} V(t^n_{k+1}, X^n_{t^n_k \wedge \cdot}) \mathbf{1}_{ \big\{(t^n_k, t^n_{k+1}] \big\}}.
		\end{align*}
		Notice that, for all fixed $\delta > 0$ and $\om \in \Om$, there exists only a finite number of $\tau^n_{k+1}$ not in $\pi_n$.
		Further, by \eqref{eq:bound_V}, the terms $\nabla_{\xr} V(\tau^n_{k+1}, X^n_{\tau^n_k \wedge \cdot})$ and $\nabla_{\xr} V(\tau^n_{k+1}, X^{n, \tau^n_{k+1}-})$ are uniformly bounded for every fixed $\om \in \Om$.
		Then by the continuity of $\nabla_{\xr} V$ and the fact that $X$ has c\`adl\`ag paths, it is easy to see that,  for every fixed $\om \in \Om$,
		\begin{equation} \label{eq:I_d2I}
			\lim_{n \to \infty} \big| I^{\delta}_n (\om)  - I_n (\om) \big| ~=~ 0. 
		\end{equation}

		{\rm (b)} Let us now assume that, for some $\delta> 0$,  $\P$ belongs to the  collection of probability measures 
		$$
			\Pc_{\delta} 
			~:=~
			\big\{ 
				\P' \in \Pc ~: \P' \big[ |\Delta X_r| \in \{0\} \cup [\delta, \infty), ~\forall r \in [0,T]  \big] = 1
			\big\},
		$$
		i.e. $X$ has only big jumps (with jump size bigger than $\delta$) under $\P$.
		As $\nabla_{\xr} V \in \Cb(\Theta)$, then
		$$
			H^n_r
			\longrightarrow 
			  \nabla_{\xr} V^-(r, X),
			~\mbox{for all}~ {r \in [s, t]},
			~\P \mbox{--a.s. for each}~\P \in \Pc_{\delta},
		$$
		in which $(\nabla_{\xr} V^-(r, X))_{r\ge 0} := (\nabla_{\xr} V(r, X^{r-}))_{r\ge 0}$ is $\F$-predictable. 
		Further, notice that one can localize the sequence of processes $( {X_{0}\1_{\{0\}}}+\sum_{k=0}^{n-1} X_{t^n_k} \mathbf{1}_{(t^n_k, t^n_{k+1}]})_{n \ge 1}$ uniformly  by using the sequence of $\F$--stopping times $\tau_m := \inf\{ t ~: |X_t| \ge m\}$, $m \ge 1$. 	Then,  by \eqref{eq:bound_V},
		the sequence $(|H^{n}|)_{n\ge 1}$ can be uniformly bounded by a locally bounded predictable process.
		By Jacod and Shiryaev \cite[Theorem I.4.31]{JacodShi}, and after possibly passing to a subsequence, 
		it follows that
		\begin{equation*} \label{eq:OptDec_big_jump}
			\int_s^t H^n_r \cdot dX_r\to  \int_s^t  \nabla_{\xr} V^-(r,X) \cdot dX_r\;\mbox{ and }\; V(t, X^n)\to V(t,X),
			~\P \mbox{--a.s.} 
		\end{equation*}  
		Therefore, \eqref{eq:OptDecomp_Xn} and \eqref{eq:I_d2I} imply that,
		for all $\P \in \Pc_{\delta}$,
		\begin{align}\label{eq:OptDec_big_jump}
			V(t, X) - V(s,X) 
			~\le~
			\int_s^t  \nabla_{\xr} V^-(r,X) \cdot dX_r,
			~~\P \mbox{--a.s.}
		\end{align} 

		{\rm (c)} We now consider $\P \in \Pc$, under which $X$ is a general semimartingale taking value in the interior of the set $E$.
		Under $\P$, $X$ can be uniquely decomposed as the sum of a continuous martingale $X^c$ and a purely discontinuous semimartingale $X^d$.
		Recall that every purely discontinuous semimartingale can be approximated uniformly, on $[0,T]$, by processes with finite variation
		({see e.g. \cite[Section I.4 and Theorem II.2.34]{JacodShi}}).
		Namely, by keeping only the (compensated) small jumps in $X^d$,
		one can find a sequence $(Z^n)_{n \ge 1}$ of purely discontinuous semimartingales, 
		together with a sequence of positive real numbers $(\delta_n)_{n \ge 1}$, 
		such that
		\begin{align} \label{eq: conv Yn to 0}
			\P \big[ |\Delta Z^n_t | < \delta_n, ~t \in [0,T] \big] = 1,
			&~\delta_n \longrightarrow 0,
			~\mbox{and}~
			\| Z^n \| +[Z^{n}]_{T} \longrightarrow 0, 
			~\P\mbox{--a.s.}, 
		\end{align}
		and $Y^n := X - Z^n$ has only jumps bigger than $\delta_n$.
		Notice that $Y^n$ may not take value in $E$ when $E \neq \R^d$.
		Let us define
		$$
			\tau_n ~:=~ \inf\{ r \ge s ~: Y^n_{r} \notin E \},
			~~
			\Yb^n_r := Y^n_r \1_{\{r < \tau_n\}} + Y^n_{\tau_n-} \1_{\{ r \ge \tau_n \}},~ r \in [s, t],
		$$
		so that $\P \circ (\Yb^n)^{-1} \in \Pc_{\delta}$.
		Then, applying \eqref{eq:OptDec_big_jump} to $\Yb^n$ {leads to}
		\begin{equation} \label{eq:Vts}
			V(t, \Yb^n) 
			~\le~
			V(s, \Yb^n) + \int_s^t  \nabla_{\xr} V^-(r, \Yb^n) \cdot d \Yb^n_r, 
			~\P \mbox{--a.s.},
			~n\ge 1.
		\end{equation}
		As $X$ takes values in  the interior of $E$, and $\| X - Y^n \| \to 0$, $\P$--a.s.,
		then, for $\P$--a.e. $\om$, there exists $n_0(\om)$ such that $\tau_n(\om) = \infty$ for all $n \ge n_0(\om)$.
		Moreover, {since} $Z^n$ is a purely discontinuous semimartingale with jumps no bigger than $\delta_n$,
		 one can localise the process, so that both $\nabla_{\xr} V^-(\cdot, \Yb^n)$, $Z^n$ and $[Z^n]_T$ are uniformly bounded.
		Taking the limit $n \to \infty$, we deduce from \eqref{eq: conv Yn to 0}, \eqref{eq:Vts} and  \cite[Theorem I.4.31]{JacodShi} that  
		\eqref{eq:OptDec_big_jump} holds true for all $\P \in \Pc$ under which $X$ {takes values}  in the interior of $E$.
	
		\vspace{0.5em}
		
		{\rm (d)} We finally consider an arbitrary $\P \in \Pc$ under which $X$ a semimartingale taking values in $E$.
		By \eqref{eq:Support_phi}, there exists a vector $e \in \R^d$ such that, for all $\eps > 0$, $X^{\eps} := X + \eps e$ is a semimartingale taking values in the interior of $E$.
		Applying  \eqref{eq:OptDec_big_jump} to $X^{\eps}$ and then letting $\eps \to 0$, 
		it follows that \eqref{eq:OptDec_big_jump} holds true for all $\P \in \Pc$.
		By the arbitrariness of $s \le t$ and $\P \in \Pc$, this proves \eqref{eq:OptDecomp} under the additional condition that $V \in \Cb^{0,1}(\Theta)$.
		
		\vspace{0.5em}

		\noindent $\mathrm{(ii)}$ 
		We can now  consider the general case without the additional condition $V \in \Cb^{0,1}(\Theta)$.
		Let $y\in \R^{d}\mapsto \phi^{\eps}(y) := \eps^{-d} \phi(\eps^{-1}y)$, $\eps>0$,  
		where $\phi$ is the smooth density function satisfying \eqref{eq:Support_phi},
		and define $V^{\eps}: \Theta \to \R$ by 
		$$
			V^{\eps}(r, \xr) 
			~:=~
			\int_{\R^d} V \big( r, \xr \oplus_r y' \big) \phi^{\eps}(y'-\xr_{r}) dy'.
		$$
		By Stokes formula and a change of variables,  
		$$
			\nabla_{\omega} V^{\eps}(r, \xr)
			~:=~
			\int_{\R^d} \eps^{-1} \Big[ V \big(r, \xr \oplus_r (\om_r + \eps y) \big)-V(r, \xr) \Big] \big( -\nabla\phi(y) \big) dy,
		$$
		 in which $\nabla \phi$ is the gradient of $\phi$.  Then,  
		$$
			V^{\eps} \in \Cb^{0,1}(\Theta),
			~~
			 V^{\eps}(r, \xr) \to V(r, \xr),
			~~\mbox{and}~~
			\nabla_{\xr} V^{\eps}(r, \xr) \to H(r, \xr),
			~~\mbox{as}~
			\eps \longrightarrow 0,\;\forall\;(r,\xr)\in \Theta, 
		$$
		where
		$$
			H(r, \xr) := \int_{\R^{d}}  \partial^{+}V(r, \xr; y) \big( -\nabla \phi(y) \big) dy,
		$$
		with 
		$$
			\partial^{+} V(r,\xr; y)
			~:=~
			\lim_{\eps\searrow 0} \frac{ V \big(r, \xr \oplus_r (\om_r+\eps y) \big) - V(t, \xr)}{\eps},~\mbox{for all}~y~\mbox{in the support of}~\phi,
		$$
		is well-defined since $V$ is Dupire-concave (see also Remark \ref{rem:DVplus} below).
		Using \eqref{eq:bound_V},
		up to a  localisation argument,
		one can assume w.l.o.g. that $(\nabla_{\xr} V^{\eps}(\cdot,X^{\cdot - }))_{\eps>0}$ is uniformly bounded.
		Then, using the decomposition result \eqref{eq:OptDec_big_jump} on $V^{\eps}$, and letting $\eps \to 0$,
		we can apply  \cite[Theorem I.4.31]{JacodShi} to conclude that, for all $\P \in \Pc$,
		\begin{align*}\label{eq:OptDec_big_jump}
			V(t, X)\le V(s,X)+\int_s^t   H(r,X^{r-}) \cdot dX_r,
			~\P \mbox{--a.s.}
		\end{align*} 
		As $\P \in \Pc$ and $s<t$  are arbitrary and $H$ does not depend on $\P$ and $s<t$, this proves the decomposition result \eqref{eq:OptDecomp}.
	
		\vspace{0.5em}
		 
		\noindent $\mathrm{(iii)}$ Finally, recalling the definition of the  super-differential of a concave function, 
		the fact that $H(r,\om)  \in \partial V(r, \om)$ is an immediate consequence of the fact that the sequence of Dupire-concave functionals $(V^{\eps})_{\eps>0}$ converges pointwise to  $V$. 
		Moreover, it is clear that the process $(H(s, X^{s-})_{s \in [0,T]}$ is $\F$--predictable,
		and is a locally bounded process by \eqref{eq:bound_V}.
	\end{proof}

	\begin{remark} \label{rem:DVplus}
		One can check that\footnote{We would like to thank Pierre Cardaliaguet who pointed out to us this identity and its proof.}
		$$
			\partial^{+} V(r,\xr; y)=\min\{y\cdot z : z\in \partial V(r,\xr)\}.
		$$
		Indeed, first, it is clear from the definition of  $\partial^{+} V(r,\xr; y)$ 
		that 
		$
			\partial^{+} V(r,\xr; y) \le \min\{y\cdot z : z\in \partial V(r,\xr)\}.
		$
		Next, let us consider $z^{\eps}\in {\rm arg}\min\{y\cdot z : z\in \partial V \big(r,\xr\oplus_{r} (\om_r+\eps y) \big)\}$, so that $z^{\eps}\cdot (\eps y) \le V \big(r,\xr\oplus_{r} (\om_r + \eps y) \big)-V(r,\xr)$. 
		By \eqref{eq:bound_V},  one can then find a sequence $(\eps_{n})_{n\ge 1}$ converging to $0$ such that $z^{\eps_{n}}\to z \in \partial V(r,\xr)$, and  $z\cdot y \le\partial^{+} V(r,\xr; y)$.
	\end{remark} 

\section{Super-hedging duality}
\label{sec:Duality}

	Let us now turn to the main motivation of this paper. 
	From Theorem \ref{thm:OptionalDecomp}, we derive in this section  a robust super-hedging problem and provide a duality result. We first state it under general abstract conditions, Theorem \ref{thm:surrep duality}, and then discuss a typical example of applications in Proposition \ref{prop: super hedging Mc+}. 
	
\subsection{Abstract framework}

	Let $\Phi : \Om \to \R$ be a payoff function and let $\Mc_0 = (\Mc(0, x))_{x \in E}$ be a family of collections of 
	probability measures $\Q$ on $\Om$ such that $X$ is a $\Q$-local martingale with $X_{0}=x$, $\Q$--a.s.
	We assume that, for all $x \in E$ and $\Q\in \Mc(0, x)$,
	\begin{equation} \label{eq:Integrability_Phi}
		\E^{\Q} \big[ \big| \Phi(X) \big| \big]  < \infty,
		~\big(\E^{\Q}[\Phi(X)^{-}|\Fc_{t}]\big)_{t\le T}
		~\mbox{is a $\Q$-martingale, and}~
		\!\! \sup_{\Q \in \Mc(0,x)} \E^{\Q} \big[ \Phi(X) \big] < \infty.
	\end{equation}

	The super-hedging price of a derivative option with payoff $\Phi(X)$ is defined by
	$$
		{\rm v}(0, x)
		~:=~
		\inf \big\{v\in \R : \exists\; H\in \Ac \mbox{ s.t. } Y_{T}^{v,H}\ge  \Phi(X),\; \Mc(0, x)-{\rm q.s.} \big\},
	$$
	in which 
	$$
		Y^{v,H}:=v+\int_{0}^{\cdot} H_{r}\cdot dX_{r}
	$$
	and $\Ac$ is the collection of all locally bounded $\F$-predictable processes such that $Y^{v,H}$ is $\Q$-a.s.~bounded from below by a $\Q$-martingale, for all $\Q\in \Mc(0, x)$.

	\vspace{0.5em}

	The aim of this section is to prove the following super-hedging duality: 
	\begin{equation} \label{eq:superhedging}
		{\rm v}(0,x) ~=~ V(0,x) ~:= \sup_{\Q\in \Mc(0, x)}\E^{\Q}[\Phi(X)].
	\end{equation}

	As usual, one can easily obtain the weak duality
	\begin{equation} \label{eq:weak_duality}
		{\rm v}(0, x) ~\ge~ V(0,x) ~:= \sup_{\Q\in \Mc(0, x)}\E^{\Q}[\Phi(X)].
	\end{equation}
	Indeed, for all $(v,H)\in \R\x \Ac$ such that $Y_{T}^{v,H}\ge  \Phi(X)$, $ \Mc(0,x)$-{\rm q.s.},
	one has $v\ge  \E^{\Q} \big[\Phi(X) \big]$ for all $\Q \in \Mc(0,x)$,
	 since $Y^{v,H}$ is a $\Q$-local-martingale bounded from below by a $\Q$-martingale, and therefore a $\Q$-supermartingale, for any $\Q\in \Mc(0,x)$,  whenever $H\in \Ac$. 

	\vspace{0.5em}

	To prove the converse inequality, we need to assume a structure condition on $\Mc(0,x)$. 
	It can be  compared to the conditions used in Biagini, Bouchard, Kardaras and Nutz \cite{BBKN13} and Nutz \cite{nutz2015OptionalDecompo} to ensure that the value function defined by the right-hand side of \eqref{eq:superhedging}, see also \eqref{eq: Vtom} below, is non-anticipative and satisfies the dynamic programming principle. 
	
	\begin{assumption} \label{assum:main} 
		There exists a family $(\Mc(t, \om))_{(t, \om) \in [0,T] \x \Om}$ of collections of probability measures on $\Om$ such that, for all $(t, \om) \in \Theta$ :
	\begin{enumerate}[{\rm (i)}]
		\item  $\Mc(0,\om)=\Mc(0, \om_0)$ and $\Mc(t, \om) = \Mc(t, \om_{t \wedge \cdot})$. 
		
		\item  For all $\Q \in \Mc(t, \om)$, $X$ is a $\Q$--local martingale on $[t,T]$ and $\Q[ X_s = \om_s, ~s \le t] = 1$.

		\item Given $\Q \in \Mc(t, \om)$ and a  $\F$--stopping time $\tau$ taking values in $[t,T]$:
		\begin{itemize}
			\item[(a)] There exists  a family  $(\Q_{\om})_{\om \in \Om}$ of r.c.p.d. of $\Q$ knowing $\Fc_{\tau}$ such that
			$$
				\Q_{\om} \in \Mc(\tau(\om), \om), 
				~\mbox{for}~
				\Q \mbox{--a.e.}~
				\om \in \Om.
			$$ 
			
			\item[(b)]
			 For all $\eps > 0$, there exists $\Q^{\eps} \in \Mc(t, \om)$ such that $\Q|_{\Fc_{\tau}} = \Q^{\eps}|_{\Fc_{\tau}}$ and
			 a family $(\Q^{\eps}_{\om})_{\om \in \Om}$ of r.c.p.d. of $\Q^{\eps}$ knowing $\Fc_{\tau}$ such that
			$$
				\E^{\Q^{\eps}} \big[\Phi(X) \big] \ge V(t, \om) - \eps
				~\mbox{and}~
				\Q^{\eps}_{\om} \in \Mc(\tau(\om), \om),
				~\mbox{for}~
				\Q \mbox{--a.e.}~ \om \in \Om,
			$$
			where
			\begin{equation}\label{eq: Vtom}
				V(t, \om) := \sup_{\Q \in \Mc(t, \om)} \E^{\Q} \big[\Phi(X) \big].
			\end{equation}
		\end{itemize}	
	\end{enumerate}
	\end{assumption}

	Assuming that the function $V$ defined above satisfies the conditions of Theorem \ref{thm:OptionalDecomp}, we can then simply apply it to deduce that there exists a $\F$-predictable process $H$,
	such that
	$Y^{V(0,x),H} \ge V(\cdot,X)$ on $[0,T]$ and in particular that $Y_{T}^{V(0, x),H}\ge  \Phi(X)$, $\Mc(0, x)$-{\rm q.s.} 
	Since, for all $\Q\in \Mc(0, x)$, $V(t, X) \ge \E^{\Q}[\Phi(X)|\Fc_{t}] \ge \E^{\Q}[\Phi(X)^{-}|\Fc_{t}]$ $\Q$-a.s.,
	and $\big(\E^{\Q}[\Phi(X)^{-}|\Fc_{t}]\big)_{t \in[0, T]}$ is a $\Q$-martingale,
	it follows that $H\in \Ac$,  and therefore that $V(0, x) \ge {\rm v}(0, x)$. 
	Together with the weak duality \eqref{eq:weak_duality}, this implies the duality result \eqref{eq:superhedging}.

	\vspace{2mm}

	To ensure that the conditions of Theorem \ref{thm:OptionalDecomp} hold, 
	we impose the following additional assumptions.
	Let {{\rm conv}$(A)$} denote the convex envelope of a set $A \subset \R^d$,
	and, for all $(t, \om) \in \Theta$ and $\eta > 0$,
	let $\delta_{\omega}$ denotes the Dirac measure at $\om$ and set
	$$
		B_{\eta} (t,\om)
		~:=~
		\big\{
			(t', \om') \in \Theta ~: t' \ge t, ~|t' - t| + \| \om_{t \wedge \cdot} - \om'_{t' \wedge \cdot} \| \le \eta
		\big\}.
	$$

	\begin{assumption} \label{assum:main2} 
		The following holds for all $(t,\om)\in [0,T) \x \Om$ and $x^{1}, x^{2}\in E$ such that $\om_{t}\in {\rm conv}\{x^{1},x^{2}\}$:
	\begin{enumerate}[{\rm (i)}] 
		\item $\delta_{\om_{t \wedge \cdot}} \in \Mc(t, \om)$,
	 
		\item for all $\eps> 0$,
		\begin{enumerate}
			\item[(a)] for all $\Q \in \Mc(t,\om)$, we can find $\eta > 0$ such that 
			for all $(t', \om') \in B_{\eta}(t,\om)$ 
			there exists $\Q' \in \Mc(t', \om')$ satisfying
				$\E^{\Q} [\Phi(X)] \le \E^{\Q'}[\Phi(X)] + \eps$.
				
			\item[(b)] there is $\eta > 0$ such that for all $(t', \om') \in B_{\eta}(t,\om)$ and   $\Q' \in \Mc(t',\om')$ we can find $\Q \in \Mc(t, \om)$ satisfying
				$\E^{\Q'} [\Phi(X)] \le \E^{\Q}[\Phi(X)] + \eps$.
		\end{enumerate}
	 
		\item
		there exists a family $(\Q_h, A^1_h, A^2_h)_{h \in (0, h_0)}$ for some $h_0 \le T-t$,
		such that $\Q_h \in \Mc(t, \om)$, $\Q_h[ X \in A^1_h \cup A^2_h ] = 1$, and
		$$
			A^i_h \subset \{ \om' \in \Om ~: \om'_s = \om_s ~\mbox{on}~[0,t],~ |\om'_s| \le |x^1| + |x^2| ~\mbox{on}~(t, t+h), ~\om'_{t+h}= x^i \},
			~i=1,2.
		$$
		Moreover, for each $\eps > 0$, $i=1,2$, $\Q \in \Mc(t, \om \oplus_t x^i)$ and $h_1 > 0$,
		there exists $h < h_1$ such that for all $\om' \in A^i_{h}$ one can find
		 $\Q' \in \Mc(t+ h, \om')$ satisfying $\E^{\Q}[\Phi(X)] \le \E^{\Q'}[ \Phi(X)] + \eps$.

	\end{enumerate}
	\end{assumption}

	Under these conditions, we can now state our main result which is an immediate consequence of the discussion above combined with Theorem \ref{thm:OptionalDecomp}, Lemma \ref{lemm:DPP} and Lemma \ref{lemm:concavity_V} below. Namely, (i) implies that $V$ is Dupire-non-increasing in time, (ii) is used to prove that $V\in \Cb(\Theta)$, while (iii) ensures that $V$ is Dupire-concave. A particular case of application will be studied in Section \ref{subsec:Example} below.

 	\begin{theorem} \label{thm:surrep duality} 
		Let Assumptions \ref{assum:main} and \ref{assum:main2} hold true.
		Assume in addition that $V$ satisfies \eqref{eq:bound_V},
		and that $y\in E \mapsto \Phi(\om \oplus_T y)$ is concave for all $\om \in \Om$.
		Then the duality \eqref{eq:superhedging} holds true,  
		and there exists a $\F$-predictable process $H \in \Ac$ such that 
		$Y^{V(0, x),H}_{T} \ge \Phi(X),~ \Mc(0, x)-\mbox{q.s.}$, for all $x\in E$.
	\end{theorem}

	\begin{remark}
		The condition that $y \mapsto \Phi(\om \oplus_T y)$ is concave for all $\om \in \Om$ is not  important  as soon as the collection $\Mc(0,x)$ is rich enough.
		For a general payoff function $\Phi: \Om \to \R$, let us denote by $\widehat \Phi: \Om \to \R$ the smallest function dominating $\Phi$ and such that $y\in E \mapsto \widehat \Phi(\om \oplus_T y)$  is concave.
		In many  situations, such as in the example of Section \ref{subsec:Example},
		one can show that
		$$
			V(t, \om) = \widehat V(t, \om) := \sup_{\Q \in \Mc(t, \om)} \E^{\Q} \big[ \widehat \Phi(X) \big],
			~\mbox{for all}~(t, \om) \in [0,T) \x \Om.
		$$
		Then, one only needs to work on $\widehat \Phi$ and $\widehat V$ to obtain the duality result \eqref{eq:superhedging} for $\widehat \Phi$,
		and then to use the weak duality \eqref{eq:weak_duality} and the above identity to deduce that \eqref{eq:superhedging} holds for $\Phi$ as well.
		See the proof of Proposition \ref{prop: super hedging Mc+} below for more details.
	\end{remark}

	The rest of this section is dedicated to the proof of the two lemmas mentioned above. 

	\begin{lemma} \label{lemm:DPP}
		The value function $V$ is non-anticipative and belongs to  $\Cb(\Theta)$. For all $(t, \om) \in \Theta$ and all $\F$-stopping times $\tau$ taking value in $[t,T]$, 
		\begin{equation} \label{eq:DPP}
			V(t, \om) 
			=
			\sup_{\Q \in \Mc(t, \om)} \E^{\Q} \big[ V(\tau, X) \big].
		\end{equation}
		In particular,  
		$V(t, \om) \ge V(t+h, \om_{t \wedge \cdot})$ for all $t \in [0,T]$, $h \in [0, T-t]$, $\om \in \Om$.
	\end{lemma}
	\begin{proof}
		$\mathrm{(i)}$. First, it is clear that $V$ is non-anticipative by Assumption \ref{assum:main}.(i).

		\vspace{0.5em}
		
		$\mathrm{(ii)}$. We next prove that $V \in \Cb(\Theta)$.
		Let $(t, \om) \in \Theta$, $(t^n, \om^n)_{n \ge 1} \subset \Theta$ be a sequence such that $t_n \searrow t$ and $\|\om^n_{t^n \wedge \cdot} - \om_{t \wedge \cdot} \| \to 0$.
		By Assumption \ref{assum:main2}(ii), for any $\eps>0$ and $\Q \in \Mc(t, \om)$ such that $\E^{\Q} [\Phi(X)] \ge V(t, \om) - \eps$,
		there exists a sequence of $(\Q_n)_{n \ge 1}$ such that $\Q_n \in \Mc(t^n, \om^n)$ and $\E^{\Q_n}[ \Phi(X)] \ge \E^{\Q} [\Phi(X)] - \eps$ for $n$ large enough. 
		This implies that
		$$
			\liminf_{n \to \infty} V(t^n, \om^n) 
			~\ge~
			\liminf_{n \to \infty} \E^{\Q_n} [ \Phi(X)]
			~\ge~ 
			\E^{\Q} [\Phi(X)] - \eps
			~\ge~
			V(t, \om) - 2 \eps.
		$$
		Next, let $(\Q'_n)_{n \ge 1}$ be a sequence such that $\Q'_n \in \Mc(t^n, \om^n)$ for all $n\ge 1$ and $\lim_{n \to \infty} \E^{\Q'_n}[ \Phi(X) ] = \limsup_{n \to \infty} V(t^n, \om^n)$.
		By Assumption \ref{assum:main2}(ii) again, for all $\eps > 0$, there exists a sequence $(\Q_n)_{n \ge 1} \subset \Mc(t, \om)$ such that
		$\E^{\Q_n} [\Phi(X)] \ge \E^{\Q_n'}[ \Phi(X)] - \eps$ for $n \ge 1$ large enough. Hence, 
		$$
			\limsup_{n \to \infty} V(t^n, \om^n)
			~=~
			\lim_{n \to \infty} \E^{\Q'_n}[ \Phi(X) ]
			~\le~
			\limsup_{n \to \infty} \E^{\Q_n}[ \Phi(X) ] + \eps
			~\le~
			V(t, \om) + \eps.
		$$
		By arbitrariness of $\eps$, one concludes that $\lim_{n \to \infty} V(t^n, \om^n) = V(t, \om)$.

		\vspace{0.5em}
		
		$\mathrm{(iii)}$. Finally, the dynamic programming principle \eqref{eq:DPP} is a direct consequence of Assumption \ref{assum:main}.(iii),
		and the fact that $V(t, \om) \ge V(t+h, \om_{t \wedge \cdot})$, for all $(t,\om)\in \Theta$ and $h\in [0,T-t]$, follows from \eqref{eq:DPP} and
		Assumption \ref{assum:main2}.(i).
	\end{proof}

	\begin{lemma} \label{lemm:concavity_V}
		The value function $V$ is Dupire-concave.
	\end{lemma}
	\begin{proof}
		We first notice that $V(T, \om) = \Phi(\om)$ by  definition, so that $y \in E \mapsto V(T, \om \oplus_T y)$  is concave, for all $\om \in \Om$.
		Let us now set $(t, \om) \in [0, T) \x \Om$ and $\om^1, \om^2$ such that $\om^1_s = \om^2_s = \om_s$ for all $s \in [0,t)$ and   $\om_t = \theta \om^1_t + (1-\theta) \om^2_t$ for some $\theta \in (0,1)$.
		Set $x^1 := \om^1_t$, $x^2 := \om^2_t$. By Assumption \ref{assum:main2}.(iii),
		there exists a family $(\Q_h, A^1_h, A^2_h)_{h \in (0, T-t)}$ such that,
		for all $h \in (0, T-t)$, one has $\Q_h \in \Mc(t, \om)$, $\Q_h[ X_{t+h} = \om^1_t] = \theta$ and $\Q_h[ X_{t+h} = \om^2_t] = 1-\theta$.
		By \eqref{eq:DPP},
		\begin{eqnarray*}
			V(t, \om)  
			&\ge&
			\E^{\Q_h} \big[ V(t+h, X) \big] \\
			&=&
			 \theta \E^{\Q_h} \big[ V(t+h, X) \big| X_{t+h} = \om^1_t \big] 
			 +
			 (1 - \theta)  \E^{\Q_h} \big[ V(t+h, X) \big| X_{t+h} = \om^2_t \big] \\
			&\ge&
			 \theta \inf_{\om' \in A^1_h} V(t+h, \om')
			 +
			 (1-\theta)  \inf_{\om' \in A^2_h} V(t+h, \om').
		\end{eqnarray*}
		Fix $i \in \{1,2\}$ and $\eps>0$.  
		Let $(\om^{h,i})_{h > 0}$ be such that $\om^{h,i} \in A^i_{h}$  and
		$$
			\inf_{\om' \in A^i_{h}} V(t+h, \om')
			~\ge~
			 V(t+h, \om^{h,i}) -\eps,
		$$
		for each $h>0$.  Let  $(\Q^i_{n})_{n \ge 1}\subset \Mc(t, \om^i)$ be such that $V(t,\om^i)=\lim_{n \to \infty} \E^{\Q^i_{n}}[ \Phi(X)]$.
		Then, by Assumption \ref{assum:main2}(iii), we can find $h_n \to 0$ and
		 a sequence $(\Q'_{h_{n}})_{ n \ge 1}$ such that 
		  $\Q'_{h_n} \in \Mc(t+h_n, \om^{h_n,i})$ and 
		$  \E^{\Q^i_n}[\Phi(X)]\le  \E^{\Q'_{h_n}}[\Phi(X)] + \eps$ for all $n\ge 1$. It follows that 
		$$
			V(t,\om^i)
			~\le~
			\limsup_{n\to \infty } \E^{\Q'_{h_n}} [\Phi(X)] + \eps
			~\le~
			\limsup_{n\to \infty } V(t+ h_n, \om^{h_n,i}) + \eps.
		$$
		Combining the above implies that
		$$
			V(t, \om) ~\ge~ \theta V(t, \om^1) + (1-\theta) V(t, \om^2).
		$$
	\end{proof}

\subsection{Example: robust hedging with positive martingales}
\label{subsec:Example}

	We now turn to a typical example of application. 
	where we consider the one-dimensional ($d=1$ for simplicity) non-negative martingales.
	Let $E = \R_+ = [0, \infty)$, so that $\Om = D([0,T], \R_+)$ and \eqref{eq:Support_phi} holds.
	Let
	$$
		\Mc^+(t, \om) 
		:=
		\big\{
			\Q ~: \Q [ X_{t \wedge \cdot} \!=\! \om_{t \wedge \cdot} ] \!=\! 1,
			~X ~\mbox{is}~ \Q \mbox{--martingale on}~ [t, T]
		\big\}.
	$$
	Given
	$$
		M_t(\om) := \sup_{0 \le s \le t} \om_s,
		~~~
		m_t(\om) := \inf_{0 \le s \le t} \om_s,
		~~
		A_t(\om) := \int_0^t \om_s \mu(d s),  \;t\le T,
	$$
	 in which $\mu$ is a finite signed measure on $[0,T]$ without atom,
	and a uniformly continuous function  $\phi: \R^4 \to \R$, we define 
	$$
		\Phi(\om) ~:=~ \phi \big(M_{T}(\om), m_{T}(\om), A_T (\om),  \om_T \big).
	$$
	We assume that there exist $K > 0$ and $\eps > 0$, such that
	\begin{equation} \label{eq:Growth_Phi}
		\big| \Phi(\om) \big| ~\le~ K \Big( 1 + \om_T + \int_0^T \om_t |\mu|(dt) \Big),
		~~\mbox{for all}~\om \in \Om,
	\end{equation}
	and, for all $0 \le M_0 \le M_1$, $0 \le m_1 \le w_1 \wedge \eps$ and $a_0, a_1\in \R$,
	\begin{equation} \label{eq:Lip_Phi}
		\Big| \phi(M_1, m_1, a_1, w_1) - \phi(M_0, 0, a_0, 0) \Big|
		\le K \big( |a_1 - a_0| + w_1 \big).
	\end{equation}
	Let us then introduce the value function
	$$
		V^+(t,  \om) := \sup_{\Q \in \Mc^+(t, \om)} \E^{\Q} \big[ \Phi(X) \big],
		~\mbox{for all}~(t, \om) \in [0,T) \x \Om,
	$$
	and $V^{+}(0,x):=V^{+}(0,x\1_{[0,T]})$ as well as $\Mc^+(0, x):=\Mc^{+}(0,x\1_{[0,T]})$, for each $x\in \R_+$.
	
	\begin{remark} \label{rem:BoundedV_DV}
		The technical conditions in \eqref{eq:Growth_Phi} and \eqref{eq:Lip_Phi} will be used to ensure Condition \eqref{eq:bound_V}.
		Namely, the growth condition in \eqref{eq:Growth_Phi} implies that 
		\begin{equation} \label{eq:Vplus_bounded}
			| V^+(t, \om) | ~\le~ K \Big( 1 + \om_t + \int_0^T \om_{t \wedge \cdot} |\mu|(dt) \Big),
		\end{equation}
		as $X$ is a non-negative martingale on $[t,T]$ under each $\Q \in \Mc^+(t, \om)$.
		Further, notice that, whenever $V^+$ is Dupire-concave, one has 
		\begin{equation} \label{eq:BoundGradient}
			\max_{z \in \partial V^+(t,\om)} |z| \le \max \big\{ |V^+(t, \om \oplus_{t} (\om_t + 1)) - V^+(t, \om)| \vee |z_0| ~: z_0 \in \partial V^+(t, \om \oplus_t 0) \big\}.
		\end{equation}
		Condition \eqref{eq:Lip_Phi} will be used to obtain a bound on $|z_0|$ for $z_0 \in \partial V^+(t, \om \oplus_t 0)$, which  then ensures that the super-gradient in $\partial V^+(t, \om)$ are also locally bounded.
	\end{remark}
	
	\begin{example}
		 Let $\phi(M, m, a, w) = f(a) + (w - g(M))_+$ for some Lipschitz function $f: \R \to \R$ and a uniformly continuous non-negative function $g: \R_+ \to \R_+$,
		then \eqref{eq:Growth_Phi} and \eqref{eq:Lip_Phi} hold true.
	\end{example}

	\begin{proposition}\label{prop: super hedging Mc+}
		Let the conditions of this subsection hold. Then, for all $x \in \R_+$,
		\begin{equation} \label{eq:duality_Mplus}
			V^{+}(0,x)
			~=~
			\inf \big\{ v \in \R ~: \exists H \in \Ac ~\mbox{s.t.}~ Y^{v, H}_T \ge \Phi(X), ~\Mc^+(0, x)-{\rm q.s.} \big\}.
		\end{equation}
		 Moreover, there exists a $\F$-predictable process $H \in \Ac$ such that 
		$$
			Y^{V^+(0,x),H}_{T} \ge \Phi(X),
			~\Mc^+(0, x)-{\rm q.s.}
		$$
	\end{proposition}
	
		\begin{remark}
		A   duality result similar to \eqref{eq:duality_Mplus} has been proved in Guo, Tan and Touzi \cite[Theorem 5.3]{GuoTT}, using the discretization technique of Dolinsky and Soner \cite{DolinskySoner} together with the S-topology technique of Jakubowski \cite{Jakubowski}.
		In \cite[Theorem 5.3]{GuoTT}, the payoff function $\Phi$ is essentially assumed to be upper semi-continuous w.r.t.~the S-topology and uniformly continuous w.r.t.~the Skorokhod topology. 
		They define the super-hedging price in terms of dynamic trading strategies $H$ that are restricted to be piecewisely constant, so that the integration $\int_0^T H_t \cdot dX_t$ can be defined $\omega$ by $\omega$, and the super-hedging property $Y^{v, H}_T \ge \Phi(X)$  also holds $\omega$ by $\omega$. 
		Our super-hedging property $Y^{v, H}_T \ge \Phi(X)$ holds in a quasi-sure sense,
		but we do not require the (semi-)continuity property w.r.t.~the S-topology (note that a uniformly continuous function of $(M_T(\om), m_T(\om))$ is generally not upper semi-continuous in $\om$ under the S-topology).
		Meanwhile, we are able to prove the existence of an optimal super-hedging strategy, which can not hold in general in the setting of \cite{GuoTT}.
		Such an optimal strategy is even given by an explicit expression, recall Remark \ref{rem: H explicit}.
	\end{remark}

	\begin{proof}{\bf (of Proposition \ref{prop: super hedging Mc+}).} 
	{\rm (i)} Let $\widehat \Phi: \Om \to \R$ be the smallest function dominating $\Phi$ and such that  $y \mapsto \widehat \Phi(\om \oplus_T y)$ is concave on $\R_+$.
	We claim that
	\begin{equation} \label{eq:claim_concave}
		V^+(t,  \om) 
		=
		\widehat V^+(t,\om) 
		:= 
		\sup_{\P \in \Mc^+(t, \om)} \E^{\P} \big[ \widehat \Phi(X) \big],
		~\mbox{for all}~(t, \om) \in [0,T) \x \Om.
	\end{equation}
	Indeed, by the definition of $\widehat \Phi$, there exists a probability space $(\Om^*, \Fc^*, \P^*)$ and a measurable map $\xi: \Om \x \Om^* \to \R_+$ such that, for all $\om \in \Om$,
	$$
		\E^{\P^*} [ \xi(\om, \cdot) ] = \om_T,
		~~\mbox{and}~~
		\widehat \Phi(\om) = \E^{\P^*} [ \Phi(\om^{T-} \oplus_T \xi(\om, \cdot)) ].
	$$
	Let $(t, \om) \in [0,T) \x \Om$ and $\P \in \Mc^+(t, \om)$, we consider the product space $(\overline \Om, \overline \Fc, \overline \P) := (\Om \x \Om^*, \Fc_T \otimes \Fc^*,  \P\x \P^*)$, and define the process
	$$
		\overline X_t ~:=~ X_t  \mathbf{1}_{\{t \in [0, T)\}} + \big( X_T + \xi \big) \mathbf{1}_{\{t = T\}},\;t\le T.
	$$
	Then,  
	$$
		\Q := \overline \P \circ \overline X^{-1} \in \Mc^+(t, \om)
		~~\mbox{and}~~
		\E^{\Q} \big[ \Phi (X) \big] = \E^{\overline \P} \big[ \Phi( \overline X) \big] = \E^{\P} \big[ \widehat \Phi(X) \big].
	$$
	This implies  \eqref{eq:claim_concave}. 	We next set $\widehat V^+(T, \om) := \widehat \Phi(\om)$ 
	and claim that  the conditions of Theorem \ref{thm:surrep duality}  are satisfied for $\widehat \Phi$ and $\widehat V$,	so that 
	\begin{align*}
		\sup_{\Q \in \Mc^{+}(0, x)} \E^{\Q} \big[ \widehat \Phi(X) \big]
		&=
		\inf \big\{v\in \R : \exists\; H\in \Ac \mbox{ s.t. } Y_{T}^{v,H}\ge \widehat \Phi(X),\; \Mc^{+}(0, x)-{\rm q.s.} \big\}	\\
		&\ge \inf \big\{v\in \R : \exists\; H\in \Ac \mbox{ s.t. } Y_{T}^{v,H}\ge   \Phi(X),\; \Mc^{+}(0, x)-{\rm q.s.} \big\}, 
	\end{align*}
	and we conclude by appealing to \eqref{eq:weak_duality} applied to $(\widehat V^{+},\widehat \Phi)$, and the existence result of Theorem \ref{thm:surrep duality}.

	\vspace{0,5em}

	{\rm (ii)} It remains to check that  the conditions of Theorem \ref{thm:surrep duality}  are satisfied for $\widehat \Phi$ and $\widehat V$.
	By construction $\widehat \Phi$ is Dupire-concave and it is straightforward to check that it inherits the uniform continuity of $\Phi$, as a function  of $(M_{T}(\om),m_{T}(\om),A_{T}(\om),\om_{T})$, as well as the bound \eqref{eq:Growth_Phi} on $\Om$. 
	First, it is easy to see that \eqref{eq:Integrability_Phi} holds true.
	In the following, we check the remaining conditions in Theorem \ref{thm:surrep duality}. 

	\vspace{0.5em}
	
	$\mathrm{(a)}$ For Assumption \ref{assum:main}, it is obvious that Items $\mathrm{(i)}$ and $\mathrm{(ii)}$ hold for $\Mc^+(t,\om)$.
	To check Item $\mathrm{(iii).(a)}$, we notice that a martingale is still a martingale under the r.c.p.d.
	For Item $\mathrm{(iii).(b)}$, one can apply  measurable selection arguments, as in e.g.~El Karoui and Tan \cite{ElKarouiTan1},
	in which the essential argument is to check that the graph set $[[\Mc^+]] := \{ (t, \om, \Q)\in [0,T] \x \Om \x \Pc(\Om)~: \Q \in \Mc^+(t,\om) \}$
	is a Borel (or only analytic) subset of $[0,T] \x \Om \x \Pc(\Om)$, where $\Pc(\Om)$ denotes the space of all Borel probability measures on $\Om$.
	Indeed, $[[\Mc^+]]$ is a Borel set as it can rewritten as
	\begin{eqnarray*}
		[[\Mc^+]] 
		\!\!&=&\!\! 
		\big\{ (t, \om, \Q) \in [0,T] \x \Om \x  \Pc(\Om)  ~: \Q[ X_{s \wedge t} = \om_{s \wedge t}] = 1, 
				~\E^{\Q}[ |X_r| + |X_s|] < \infty, \\
		&&~~~~~~~~~~~~~~~~~~~~~~~~~~~~~~~~~~~~~~~~~~~~~
			\E^{\Q} \big[ (X_{t\vee s} - X_{t \vee r}) \xi \big] = 0,
			~\mbox{for all}~(r, s, \xi) \in L
		\big\},
	\end{eqnarray*}
	where $L := L_1 \cup L_2$, $L_1$ is a countable dense subset of 
	$$
		\big\{(r, s, \xi) 
		~: 
		0 \le r < s \le T, 
		~\xi ~\mbox{being bounded continuous and}~\Fc_r\mbox{--measurable} 
		\big\},
	$$
	and $L_2$ is a countable dense subset of
	$$
		\big\{(r, T, \xi) ~: 0 \le r < T,~~\xi ~\mbox{being bounded continuous and}~\Fc_r\mbox{--measurable} 
		\big\}.
	$$
	
	$\mathrm{(b)}$ For Assumption \ref{assum:main2}, Item $\mathrm{(i)}$ is clearly true.
	To check Item $\mathrm{(ii)}$, we shall use constructions that preserve the max, the min, the $T$-value and the integral w.r.t.~$\mu$ of a given path up at a uniform distance, possibly up to an event with vanishing probability. The uniform continuity of $\om\mapsto \widehat \Phi(\om)$ as a function of $(M_{T}(\om),m_{T}(\om),A_{T}(\om),\om_{T})$ will then allow us to conclude.

	\vspace{0.5em}

	Let us first consider $(t, \om) \in [0,T) \x \Om$, $\eta > 0$, and
	$(t', \om') \in B_{2\eta}(t,\om)$, i.e.~$t' \ge t$ and $|t' - t| + \| \om_{t \wedge \cdot} - \om'_{t' \wedge \cdot} \| \le 2 \eta$.
	For simplicity, we can assume that $t=0$, $t' = \eta$ and $\|\om'_{\eta \wedge \cdot} - \om_{0\wedge \cdot}\| \le \eta$.
	Let us fix $\Q \in \Mc^{+}(0, \om)$.
	When $0 \le \om_0 \le \om'_{\eta} $, let us construct a process $\overline X^{\eta}$ as follows:
	$$
		\overline X^{\eta}_s := 
		\begin{cases}
			\om'_s, &\mbox{when}~s \in [0,\eta),\\
			\om'_{\eta} + X_{2(s-\eta)} - X_0,& \mbox{when}~ s \in [\eta, 2 \eta),\\
			\om'_{\eta} + X_s - X_0,&\mbox{when}~s \in [2\eta, T].
		\end{cases}
	$$
	When $0 < \om'_{\eta} < \om_0$, let us construct the process $\overline X^{\eta}$ as follows:
	by possibly enlarging the space, we consider a random variable $\xi_{\eta}$, independent of $X$ under $\Q$, such that $\Q[\xi_{\eta} = \om_0] + \Q[ \xi_{\eta} = 0] = 1$ and $\E^{\Q} [\xi_{\eta}] = \om'_{\eta}$,
	and we set
	$$
		\overline X^{\eta}_s := 
		\begin{cases}
			\om'_{\eta \wedge s}, &\mbox{when}~s \in [0,2\eta),\\
			\xi_{\eta}, &\mbox{when}~s = 2 \eta,\\
			X_{3(s-2\eta)} \mathbf{1}_{\{\xi_{\eta} = \om_0\}},& \mbox{when}~ s \in (2\eta, 3 \eta),\\
			X_s \mathbf{1}_{\{\xi_{\eta} = \om_0\}},&\mbox{when}~s \in [3\eta, T].
		\end{cases}
	$$
	Notice that  $\Q[ \xi_{\eta} = 0]$ converge  to $0$, as $\eta \to 0$.
	Then, in both cases, $\Q'_{\eta} := \Q \circ (\overline X^{\eta})^{-1} \in \Mc^{+}(t', \om')$ 
	and $\lim_{\eta \to 0} \E^{\Q}\big[\widehat \Phi(\overline X^{\eta}) \big] = \E^{\Q} \big[ \widehat \Phi(X) \big]$ (recall that $\widehat \Phi$ is uniformly continuous).
	Thus, for all $\eps>0$, $\E^{\Q} [\widehat \Phi(X)] \le \E^{\Q'_{\eta}}[ \widehat \Phi(X)] + \eps$ for $\eta > 0$ small enough.
	This proves Item (ii).(a) of  Assumption \ref{assum:main2}.
	
	\vspace{0.5em}
	
	Next, let $(t, \om) \in [0,T) \x \Om$ be such that $\om_t \ge 0$, $\eta > 0$, $(t', \om') \in B_{2 \eta}(t,\om)$.
	W.l.o.g., let us assume that $t= 0$ and $t' = \eta$.
	Then, for each $\Q'_{\eta} \in \Mc^+(\eta, \om')$, we construct a process $\overline X^{\eta}$ as follows.
	When $\om_0 \ge \om'_{\eta}$, we set
	$$
		\overline X^{\eta}_s := 
		\begin{cases}
			\om_0, &\mbox{when}~s \in [0, \eta),\\
			\om_0 + X_s - \om'_{\eta}, &\mbox{when}~s \in [\eta, T].
		\end{cases}
	$$
	When $0 < \om_0 < \om'_{\eta}$, by possible enlarging the space, we consider a random variable $\xi_{\eta}$, independent of $X$ under $\Q'_{\eta}$, such that $\Q'_{\eta}[\xi_{\eta} = \om'_{\eta}] + \Q'_{\eta}[\xi_{\eta} = 0] = 1$ and $\E^{\Q'_{\eta}} [\xi_{\eta}] = \om_0$,
	and we set
	$$
		\overline X^{\eta}_s := 
		\begin{cases}
			\om_0, &\mbox{when}~s \in [0, \eta),\\
			X_s  \mathbf{1}_{\{\xi_{\eta} = \om'_{\eta}\}}, &\mbox{when}~s \in [\eta, T].
		\end{cases}
	$$
	In both cases, one can check that $\Q_{\eta} := \Q'_{\eta} \circ (\overline X^{\eta})^{-1} \in \Mc^{+}(t, \om)$, 
	and $\lim_{\eta \to 0} \E^{\Q'_{\eta}}[ \widehat \Phi(X)] - \E^{\Q_{\eta}}[ \widehat \Phi(X)] = 0$. 
	This shows that, for all $\eps>0$,
	$\E^{\Q'_{\eta}}[ \widehat \Phi(X)] \le \E^{\Q_{\eta}}[ \widehat \Phi(X)] + \eps$ for $\eta > 0$ small enough.
	
	\vspace{0.5em}

	When $0 = \om_0 \le \om'_{\eta} \le \eta $, there exists only one element $\Q_0 \in \Mc^+(0,\om)$ under which $X_s \equiv 0$ for $s \in [0,T]$.
	Notice that, for any sequence $(\eta_n, \Q'_n)_{n \ge 1}$ 
	such that $\eta_n \to 0$ (so that $\om'_{\eta_n} \to 0$) and $\Q'_n \in \Mc(\eta_n, \om'_{\eta_n})$,
	one can construct, on an abstract probability space $(\Om^*, \Fc^*, \P^*)$, a sequence $(X^{*,n})_{n \ge 0}$ such that $\P^* \circ (X^{*,0})^{-1} = \Q_0 \circ X^{-1}$ and $\P^* \circ (X^{*,n})^{-1} = \Q'_n \circ X^{-1}$.
	As $X^{*,n}$ is a martingale on $[\eta_n, T]$, and $X^{*,n}_{\eta_n} = \om'_{\eta_n}  \to 0$, 
	it follows from Doob's inequality that $\sup_{0 \le s \le T} |X^{*,n}_s| = \sup_{0 \le s \le T} |X^{*,n}_s - X^{*,0}_s| \to 0$ in probability.
	By the martingale property, one has $\E^{\P^*} \big[ |X^{*,n}_T| + \int_0^T |X^{*,n}_t| |\mu|(dt) \big] \to 0$,
	and hence the sequence $\big(  |X^{*,n}_T| + \int_0^T |X^{*,n}_t| |\mu|(dt) \big)_{n \ge 1}$ is uniformly integrable.
	Then, by \eqref{eq:Growth_Phi}, the sequence $(\widehat \Phi(X^{*,n}))_{n \ge 1}$ is also uniformly integrable.
	It follows that $\E^{\P^*} [ \widehat \Phi(X^{*,n}) ] \to \E^{\P^*} [ \widehat \Phi(X^{*,0})]$, or equivalently $\E^{\Q'_n}[ \widehat \Phi(X)] \to \E^{\Q_0}[ \widehat \Phi(X)]$,
	which is enough for Item (ii).(b) of  Assumption \ref{assum:main2}.

	\vspace{0.5em}
	
	$\mathrm{(c)}$ Let us then check Item $\mathrm{(iii)}$ of Assumption \ref{assum:main2}. 
	We use a similar type of construction as in step (b) above.
	Let $(t, \om) \in [0,T) \x \Om$ and $x^1, x^2 \in \R_+$ be such that $\om_t \in \mathrm{conv}\{x^1, x^2\}$.
	For each $h > 0$, $i=1,2$, we define $A^i_h$ by
	$$
		A^i_h := \big\{
			\om' \in \Om ~: \om' = \om_{t \wedge \cdot} ~\mbox{on}~[0,t+h),~\om'_{t+h} = x^i
		\big\},
	$$
	and let $\Q_h \in \Mc^{+}(t, \om)$   be such that $\Q_h[ X_s = \om_t, ~s \in [t, t+h) ] =1$ and $\Q_h[ X_{t+h} = x^1] + \Q_h[ X_{t+h} = x^2] = 1$.
	Then, for each $h > 0$ and  $i=1,2$, we define $\overline X^{h,i}$ by
	$$
		\overline X^{h, i} _s 
		~:=~ 
		\om_{t \wedge s} \mathbf{1}_{\{s \in [0,t+h)\}}
		+
		\big( x^i + X_{t + 2(s-t-h)} - X_t \big) \mathbf{1}_{\{s \in [t+h, t+2h)\}}
		+
		\big( x^i + X_{s} - X_t \big) \mathbf{1}_{\{s \in [t+2h, T] \}}.
	$$
	For every $\Q \in \Mc^+(t, \om^{t-}\oplus_t x^i)$, we notice that $\Q'_{h} := \Q \circ (\overline X^{h, i})^{-1} \in \Mc^+(t+h, \om_{t\wedge \cdot} \oplus_{t+h} x^i)$
	and that
	$\lim_{h \to 0} \E^{\Q'_{h}} [\widehat \Phi(X)] = \lim_{h \to 0} \E^{\Q} [ \widehat \Phi(\overline X^{h,i})] = \E^{\Q}[\widehat \Phi(X)]$,
	which is enough to conclude that Item $\mathrm{(iii)}$ of Assumption \ref{assum:main2} holds true.

	\vspace{0.5em}

	$\mathrm{(d)}$ Finally, we prove that $V^+$ satisfies \eqref{eq:bound_V}.
	As discussed in Remark \ref{rem:BoundedV_DV}, the growth condition \eqref{eq:Growth_Phi} implies the locally boundedness of the function ${V^+}$, c.f. \eqref{eq:Vplus_bounded}.
	Further, in view of \eqref{eq:BoundGradient}, it is enough to prove that $y \mapsto V^+(t, \om \oplus_t y)$ is Lipschitz on $[0, \eps]$  for some $\eps > 0$.
	This is true since, for $y \in [0, \eps]$,
	$$
		\big| V^+(t, \om \oplus_t y) - V^+(t, \om \oplus_t 0) \big|
		\le\!\!
		\sup_{\Q \in \Mc^+(t, \om \oplus_t y)} 
		K \E^{\Q} \Big[
		|A_T - A_t| + X_T
		\Big]
		~\le~
		2 K \big(1 \vee |\mu|([t,T]) \big) y.
	$$
	\end{proof}

\bibliographystyle{plain}

 \end{document}